\documentclass[11pt]{article}

\usepackage{amsmath,amssymb,amsthm,epsfig}

\newcommand{\circk}{\mathop{K}\limits^\circ}
\newcommand{\vp}{\varepsilon}
\newcommand{\n}{\noindent}

\newtheorem{thm}{Theorem}
\newtheorem{pro}[thm]{Proposition}

\theoremstyle{remark}
\newtheorem{rem}{Remark}

\begin{document}

\begin{center}
{\bf LEVI FOLIATIONS IN PSEUDOCONVEX BOUNDARIES AND VECTOR FIELDS THAT
COMMUTE APPROXIMATELY WITH $\pmb{\bar\partial}$}
\footnote{Research supported in part by
NSF grant DMS-9801539}
\end{center}
\bigskip

\begin{center}
Emil J.~Straube and Marcel K. Sucheston\footnote{Marcel K.~Sucheston died
tragically on April 24, 2000. The surviving author dedicates this paper to his
memory.}
 \end{center} \vspace{.5in}

\n  Abstract: \ Boas and Straube ([\ref{B-S1}]) proved a general sufficient
condition for global regularity of the $\bar\partial$-Neumann problem in terms
of families of vector fields that commute approximately with $\bar\partial$.
In this paper, we study the existence of these vector fields on a compact
subset in the boundary whose interior is foliated by complex manifolds. This
question turns out to be closely related to properties of interest from the
point of view of foliation theory.
\bigskip

\centerline{1. INTRODUCTION}

In [\ref{B-S1}], Boas and Straube established a general sufficient condition
for global regularity of the $\bar\partial$-Neumann problem in terms of
(families of) vector fields that commute approximately with $\bar\partial$.
They then showed that when the infinite type points lie in certain
submanifolds, in particular in complex submanifolds, of the boundary, then
there is an element in the first De~Rahm cohomology of the submanifold that is
the obstruction to the existence of the required vector fields: \ these vector
fields exist if and only if this cohomology class is zero ([\ref{B-S1}], Main
Theorem and Remark~5). In the case of a complex submanifold of the boundary,
this cohomology class has a geometric interpretation as a measure of the
winding of the boundary around the submanifold ([\ref{B-F1}]). 

In this paper,
we study the situation where there is a Levi foliation in the boundary. More
precisely, we assume that the set $K$ of infinite type points (which is
compact) is the closure of its interior $\circk$ (in the boundary), and that
all these points are Levi flat (the Levi form vanishes identically). By the
Frobenius theorem, $\circk$ is foliated by complex manifolds of codimension 1;
this foliation is usually referred to as the Levi foliation of $\circk$.
We start in section~2 by briefly recalling
from [\ref{B-S1}] what is needed concerning vector fields that commute
approximately with $\bar\partial$. In section~3, we assume that the domain
$\Omega$ is in $\mathbb{C}^2$, and that the boundary of $K$ (relative to
$b\Omega$) is a smooth surface with only isolated generic complex tangencies.
Under suitable hypotheses, we show that the (families of) vector fields from
[\ref{B-S1}] exist (Theorem~1). The cohomology class mentioned above, on each
leaf of the Levi foliation, turns out to coincide with the infinitesimal
holonomy of the leaf.  Thus the question of existence of vector fields
commuting approximately with $\bar\partial$ is closely related to foliation
theoretic properties of the Levi foliation. We elaborate on this theme in
section~4; in particular, Sullivan's theory of foliation currents ([\ref{Su1}],
[\ref{Su2}], [\ref{C-C}]) yields a cohomological necessary and sufficient
condition for the existence of the vector fields (Theorem~\ref{thm3}).

\bigskip

\centerline{2. VECTOR FIELDS THAT COMMUTE APPROXIMATELY WITH
$\bar\partial$} \label{sec2}
\medskip

On a bounded pseudoconvex domain, the $\bar\partial$-Neumann operator $N_q$
$(1\le q\le n)$ is the inverse of the complex Laplacian $\bar\partial
\bar\partial^*+ \bar\partial^*\bar\partial$ acting on $(0,q)$-forms. For
background on the $\bar\partial$-Neumann problem, we refer the reader to
[\ref{F-K}] and to the recent survey [\ref{B-S2}]. We say that the
$\bar\partial$-Neumann problem is globally regular in Sobolev spaces on a
domain $\Omega$ if the operators $N_q$ map the spaces $W^s_{0,q}(\Omega)$ of
$(0,q)$-forms with coefficients in $W^s(\Omega)$ continuously to themselves
for $s\ge 0$, $1\le q\le n$. Here, $W^s(\Omega)$ denotes the usual
$\mathcal{L}^2$-Sobolev space of order $s$.

In [\ref{B-S1}], Boas and Straube showed that a sufficient condition for global
regularity is the existence of families of vector fields that have certain
approximate commutation properties with $\bar\partial$. Let $\Omega$ be a
smooth $(C^\infty)$ bounded pseudoconvex domain in $\mathbb{C}^n$ with
defining function $\rho$. Suppose there is a constant $C>0$ such that for
every $\vp>0$ there exists a vector field $X_\vp$ of type (1,0) whose
coefficients are smooth in a neighborhood $U_\vp$ in $\mathbb{C}^n$ of the
set of boundary points of $\Omega$ of infinite type (which is compact,
[\ref{DA}])  and such that
\begin{equation}\label{eq1}
|\arg X_\vp\rho| <\vp\quad \text{and}\quad C^{-1} < |X_\vp\rho|<C \text{ on }
U_\vp
\end{equation}
and
\begin{equation}\label{eq2}
|\partial\rho([X_\vp, \partial/\partial \bar z_j])| < \vp \text{ on } U_\vp,
\qquad 1\le j\le n.
\end{equation}
(\ref{eq2}) says that the normal (1,0)-component of the commutators
$[X_\vp,\partial/\partial\bar z_j]$ should have modulus less than $\vp$, i.e.\
$X_\vp$ commutes approximately with $\partial/\partial \bar z_j$, $1\le j\le
n$, modulo terms of type (0,1) and terms that are complex tangential. These
terms are benign (they can be absorbed) as far as Sobolev estimates for the
$\bar\partial$-Neumann problem are concerned (see [\ref{B-S3}], p.~83, formulas
(\ref{eq2}) and (\ref{eq3}) for a precise statement), and so we still say
that the fields $X_\vp$ commute approximately with $\bar\partial$, although
only the normal (1,0)-component of the commutator is required to be small. We
also note that in (\ref{eq2}), it suffices to consider commutators with
tangential fields of type (0,1), as one can always extend fields from to
boundary of $\Omega$ to the inside suitably (see [\ref{B-S3}], p.~86).

When computing the normal (1,0)-components of the commutators in (\ref{eq2}),
a certain 1-form appears naturally. Let $\eta :=
i(\partial\rho-\bar\partial\rho)$, and $T := -i(L_n-\bar L_n)$, where $L_n =
(2/|\nabla\rho|^2) \sum\limits^n_{j=1} \frac{\partial\rho}{\partial\bar z_j}
\frac\partial{\partial z_j}$. $L_n$ is the complex normal, normalized so that
$\eta(T) \equiv 1$ on $b\Omega$. Let $\alpha := -\mathcal{L}_T\eta$, minus the
Lie derivative of $\eta$ in the direction of $T$ ([\ref{DA}], [\ref{B-S1}];
note that $\eta$ and $T$ as defined here differ by factors of $i$ and $-i$,
respectively, from the notation in [\ref{DA}], [\ref{B-S1}], so that both
$\eta$ and $T$ are real. This does not affect $\alpha$.). If $Y$ denotes a
(local) section of $T^{0,1}(b\Omega)$, then
\begin{equation}\label{eq3}
\alpha(Y) = 2\partial\rho([L_n,Y])
\end{equation}
(see [\ref{DA}], p.~92, [\ref{B-S1}], section~2). Let $\bar L_1,\ldots, \bar
L_{n-1}$ be local sections of $T^{0,1}(b\Omega)$ which span the complex
tangent space to $b\Omega$ near a boundary point $P$. Near $P$, a vector field
$X$ on $b\Omega$ of type (1,0), with $X\rho \ne0$ on $b\Omega$, can be
written as $X = e^h L_n + \sum\limits^{n-1}_{j=1} a_jL_j$ (the $e^hL_n$ term
is actually global). We have (note that $L_n\rho = 1/2$)
\begin{align}\label{eq4}
\partial\rho([X,\bar L_k]) &= e^h\left[\frac12 \bar L_kh + \partial\rho ([
L_n, \bar L_k])\right] + \sum^{n-1}_{j=1} a_j\partial\rho([L_j,\bar
L_k])\\ 
&= \frac12 e^h[\bar L_kh + \alpha(\bar L_k)] + \sum^{n-1}_{j=1} a_j 
\partial\rho([L_j,\bar L_k]).\nonumber
\end{align}

In general, one can use both $h$ and the $a_j$'s to adjust the commutators ([\ref{B-S3}], [\ref{B-S1}], [\ref{B-S2}], Remark 1 below).
Assume now there is a Levi foliation in the boundary, as discussed in
section~1. We denote this foliation by $\mathcal{F}$. If $\bar L_k$ is tangent
to the leaves of $\mathcal{F}$, then $\partial \rho([L_j,\bar L_k]) = 0$,
$1\le j\le n-1$ (because $\bar L_k$ is a null direction of the Levi form,
which is positive semi-definite on the complex tangent space), so we only have the function $h$ to work with. If
there exists a real function $h$, smooth in a neighborhood of the set $K$ of
points of infinite type, such that \begin{equation}\label{eq5}
dh/_{\text{leaf}} = -\alpha/_{\text{leaf}}
\end{equation}
for all leaves of the Levi foliation, then the field $X := e^hL_n$ will have
properties (\ref{eq1}) and (\ref{eq2}) above in a suitably small
neighborhood $U_\vp$ of $K$ (by continuity). We note that (\ref{eq5}) requires
of course that the restriction of $\alpha$ to any leaf of the Levi foliation is
closed, but this was shown to be the case in [\ref{B-S1}] (see the Lemma in
section~2). Furthermore, the definition of $\alpha$  involves a choice of
defining function. However, it is easy to see that $\alpha$'s resulting from
different defining functions differ by an exact form on the leaves of
$\mathcal{F}$, so that whether or not (\ref{eq5}) can be solved does not
depend on the defining function. Indeed, if $\tilde\rho=e^g\rho$, then one
computes that the restriction of $\alpha$ to leaves of $\mathcal{F}$ changes by
$dg$.

The following property of $\alpha$ will be useful.  In terms
of $\eta$, the Frobenius condition reads $d\eta\wedge\eta=0$. Consequently,
there exists a one form $\beta$ such that $d\eta =\beta\wedge\eta$. In fact,
$\alpha$ is such a form, that is $d\eta=\alpha\wedge\eta$ on $K$. For this,
see for example the discussion in chapter~2 of [\ref{T}], in particular
Proposition~2.2.\medskip

\begin{rem}\label{rem1}
At this point, we can easily consider the case of a Levi foliation of higher
codimension. That is, assume that the rank of the Levi form is constant on
$K$, but possibly greater than zero (which is the flat case). If that rank is
$p$, then $\circk$ is foliated by complex manifolds of dimension $n-1-p$. If
there is a function $h$ (smooth in a neighborhood of $K$) satisfying
(\ref{eq5}), then vector fields $X_\vp$ satisfying (\ref{eq1}) and (\ref{eq2})
may be obtained using the complex tangential term $\sum\limits^{n-1}_{j=1}
a_jL_j$ to adjust the commutators when $L_k$ is tangential to $b\Omega$, but
is not tangential to a leaf of the foliation (note that in this case,
$\partial\rho([L_k, \bar L_k])\ne 0$, so that we do get a contribution in
(\ref{eq4}) from these terms). Details on how to obtain these fields locally,
and on how to patch the local fields to fields defined in neighborhoods
$U_\vp$ of $K$, may be found in [\ref{B-S3}], proof of the Lemma, and
[\ref{B-S1}], section~3.
\end{rem}
\medskip

\centerline{3. DOMAINS IN $\mathbb{C}^2$}

We assume now that $\Omega$ is in
$\mathbb{C}^2$.  Suppose that $K$ is
connected and has a smooth boundary; more precisely, $K$ is the
closure, relative to $b\Omega$, of the smoothly bounded domain (as a subset of
$b\Omega)$ $\circk$. Denote by $\Gamma$ the boundary (relative to $b\Omega$)
of $K$  (or $\circk$). Then $\Gamma$ is a smooth compact orientable surface
embedded in $\mathbb{C}^2$. Recall that a generic complex tangency of $\Gamma$
is either elliptic or hyperbolic: \ nearby, $\Gamma$ can be written, in
suitable coordinates, as $w = |z|^2 + \lambda \text{ Re } z^2 + O(|z|^3)$,
with $0 \le \lambda<1$ (the elliptic case) or $1<\lambda$ (the hyperbolic
case). This goes back to Bishop's paper [\ref{B}], as does the index formula 
$\#e-\#h=2-2g$, where $\#e$ and $\#h$ denote the number of elliptic and
hyperbolic tangencies, respectively, and $g$ is the genus of $\Gamma$.
(Compare the introduction in [\ref{F}] for information on various indices
relevant for real surfaces inside complex surfaces.) In our situation, the
following point of view is useful. The Levi foliation on $\circk$ (still
denoted by $\mathcal{F}$) is induced by the  1-form $\eta$, and there is
a corresponding 1-dimensional foliation, with singularities, on $\Gamma$,
induced by restricting $\eta$ to $\Gamma$. The singularities (i.e.\ the
points of $\Gamma$ where $\eta/\Gamma=0$) are precisely the complex
tangencies. Moreover, this foliation is orientable (since $\Gamma$ is), and so
yields a vector field with singularities at the complex tangencies. The index
is 1 at an elliptic point, and $-1$ at a hyperbolic point. Thus the above index
formula results in this case from the Poincar\'e-Hopf index theorem (compare
also [\ref{B-G}], section~3 for a discussion along these lines). 

At a hyperbolic point of $\Gamma$, there are, locally, two leaves of
$\mathcal{F}$ that meet (they can be constructed by flowing along suitable
complex tangential curves through the point). We make the (global) assumption
that these two local leaves are distinct globally.

For terminology from foliation theory, in particular for the notion of
(germinal) holonomy and of infinitesimal holonomy, respectively, of a leaf, we
refer the reader to [\ref{C-N}], [\ref{C-C}].

\begin{thm}\label{thm1}
Let $\Omega$ be a smooth bounded pseudoconvex domain in $\mathbb{C}^2$.
Suppose that the set $K$ of infinite type points of $b\Omega$ is smoothly
bounded (in $b\Omega$) and that its boundary $\Gamma$ has only isolated
generic complex tangencies. Assume that the two leaves meeting at a hyperbolic point of
$\Gamma$ have no other hyperbolic points in their closure (in $K$).
If each leaf of the Levi foliation is closed (in $\circk$) and has trivial
infinitesimal holonomy, then the vector fields discussed in section~2 exist.
Consequently, the $\bar\partial$-Neumann problem on $\Omega$ is globally
regular. 
\end{thm}

We now comment further on the assumptions.

\begin{rem}\label{rem2}
Note again that it is understood that the two  local leaves meeting at a
hyperbolic point are distinct globally. The assumption that leaves have at most
one hyperbolic point in their closure (relative to $K$) is also a kind of a
genericity assumption; it  simplifies the non-Hausdorff structure of the leaf
space $\circk/\mathcal{F}$. It will be clear from the proof of
Theorem~\ref{thm1} how to handle more general situations, provided certain
``compatibility'' conditions (arising from the fact that the local
determination of $h$ near hyperbolic points propagates along leaves, by virtue
of the relation (\ref{eq5})) are met.  \end{rem}

\begin{rem}\label{rem3}
In order to solve $dh/_{\text{leaf}} = -\alpha/_{\text{leaf}}$, we need of
course that $\alpha$ represents the trivial De~Rham cohomology class on each
leaf. However, the vanishing of this cohomology class is equivalent to the
triviality of the infinitesimal holonomy of the leaf. In fact, the
infinitesimal holonomy of a leaf $L$ can be viewed naturally as an element
in the first De~Rham cohomology of $L$, and then it {\em is\/} the class of
$\alpha$ ([\ref{C-C}], Example~2.3.15). In this context, compare also
Corollary~2 in the appendix of [\ref{Ba1}]; the computations in section~2 above
(see in particular (\ref{eq4})) show that if $e^hL_n$ is holomorphic on $L$,
$dh/_L + \alpha/_L = 0$ (taking into account that both $h$ and $\alpha$ are
real). It will be seen in section~4 that the existence of a function $h$ as
above is equivalent to the existence of a (smooth) {\em closed\/} one form
that defines the Levi foliation and thus we get that with the additional
assumptions in Theorem~\ref{thm1}, the holonomy itself of each leaf is trivial
(see e.g. [\ref{T}], Theorem~3.29). For completeness, we point out that the
condition of vanishing holonomy (hence vanishing infinitesimal holonomy) is
trivially satisfied if the leaves are simply connected (that is, they are
analytic discs, by the classification of simply connected Riemann surfaces).
\end{rem}

\begin{rem}\label{rem4}
An example where there is a leaf with nontrivial holonomy may be found in
[\ref{B-F1}], p.~21. Consider 
$$S_R := \{(z,w) \in \mathbb{C}^2 \mid |z|\le 1, 1\le |w|\le R, \text{
Re}(ze^{i\log|w|^2}) = 0\}.$$
This is a Levi flat hypersurface in $\mathbb{C}^2$. It is shown in [\ref{B-F1}]
how to ``cap off'' $S_R$ to obtain a pseudoconvex domain $\Omega\subseteq
\mathbb{C}^2$ such that $S_R\subseteq b\Omega$ and $b\Omega\backslash S_R$ is
strongly pseudoconvex.  $\Omega$ is a worm domain, but with the critical
annulus fattened. The leaves of the Levi foliation on $S_R$ are given by
$$L(c) := \{z=icw^{-2i}\mid 1< |w|<R, |z|<1\},$$
where $c\in\mathbb{R}$. On $L(0),\alpha$ does not represent the trivial
cohomology class; its integral around a concentric circle inside $L(0)$ (which
is an annulus) is nonzero. Also note that the leaves $L(c)$ for $c\ne 0$ are
not closed (they spiral towards $L(0)$). Barrett's arguments carry over from
the case of the ``standard'' worm domains ([\ref{Ba-2}]) to show that the
conclusion of Theorem~\ref{thm1} fails for such $\Omega$.  \end{rem}

\begin{rem}\label{rem5}
Given the results of [\ref{B-S3}], it is natural to ask whether in the
situation of Theorem~\ref{thm1}, there is a defining function that is
plurisubharmonic on the boundary. It is not hard to see that $e^h\rho$ ($h$ in
section~2) gives a defining function that is plurisubharmonic at points of $K$
if $h$ is extended from the boundary with suitable normal derivative. While
this implication does not hold in general ([\ref{B-S1}], Remark~3), in the case
at hand it is a special case of results in [\ref{S-S}].
\end{rem}
\medskip

To prove Theorem~\ref{thm1}, it suffices, by the discussion in section~2, to
find a function $h\in C^\infty(K)$ such that
\begin{equation}\label{eq6}
dh/_{\text{leaf}} = -\alpha/_{\text{leaf}.}
\end{equation}
We will construct $h$ on $\circk$ satisfying (\ref{eq6}) and in addition
having all its derivatives bounded on $\circk$. Such a function will then
extend to a function in $C^\infty(K)$ (since $\circk$ has smooth boundary).
The construction proceeds in three stages: \ in saturated neighborhoods of the
complex tangencies of $\Gamma$, in a saturated neighborhood of the closure of
the remaining points, and, finally, patching to obtain a globally defined $h$.
In the last step, the key is that we can use cutoff functions which are
constant on the leaves (i.e.\ we only need to cut off in the direction
transverse to the Levi foliation), so that (\ref{eq6}) is not affected (this
technique occurs in a related context already in [\ref{B-F1}]).

We start with a hyperbolic point $P$. There is a diffeomorphism of a
neighborhood of $P$ (in $b\Omega$) onto a neighborhood of $O\in \mathbb{R}^3$
such that in these coordinates, $K$ is given locally by $x_3\le x^2_1-x^2_2$
(by the Morse theorem). The complex tangent plane to $b\Omega$ at $P$
corresponds to the $(x_1,x_2)$-plane. Denote by $L_1$ the leaf of $\mathcal{F}$
that meets $O$ (i.e.\ $P$) from the side where $x_1<0$. (Near the origin,
this leaf is obtained by flowing along complex tangential curves emanating
from 0; alternatively, one starts at points on the (unique) complex tangential
curve in the $(x_1,x_3)$-plane through 0, and then flows along complex
tangential directions transverse to it.)  Then if $a>0$ is small enough, there
is a unique $t_0$ such that $(-a,0,t_0)\in L_1$. Then, for $t$ close enough
to $t_0$, the curve $t\mapsto (-a,0,t)$ is transversal to the leaves of
$\mathcal{F}$, and it meets each leaf at most once. This last property is
clear locally, but it also holds globally by [\ref{H}], Proposition~3.1 (note
that the assumption that the first Betti number of the underlying manifold, in
our case $\circk$, is finite, is satisfied). On these leaves, we define $h$
(globally) by the requirements that
\begin{align}\label{eq7}
h(-a,0,t) &= 0,\\
\intertext{and}
\label{eq8}
dh/_{\text{leaf}} &= -\alpha/_{\text{leaf}}.
\end{align}
We use here that
$\alpha$ represents the trivial De~Rham cohomology class on each leaf, or,
equivalently, that each leaf has trivial infinitesimal holonomy (see
Remark~\ref{rem3} above). Now consider the leaf $L_2$ which meets 0 from the side
where $x_1>0$. There is a unique $t_1$, such that $(a,0,t_1)\in L_2$ (again: \
$a>0$ small enough), and if $t$ is close enough to $t_1$, the curve $t\mapsto
(a,0,t)$ is transversal to $\mathcal{F}$ and meets each leaf at most once. If
$t<t_1$, $(a,0,t)$ belongs to a leaf that passes through a point $(-a,0,t')$
for some $t'<t_0$. On these leaves, $h$ was defined above. Note that the
function $h(a,0,t)$ is {\em smooth up to\/} $t=t_1$ (from the side
$t<t_1$). To see this, note that $h(a,0,t)$ may be obtained by line integrals
of $\alpha$ along integral curves of the vector field induced in the
$(x_1,x_3$)-plane (near the origin) by intersecting  with the complex
tangent space. Both this vector field and $\alpha$ are smooth in a {\em
full\/} neighborhood of the origin (in the $(x_1,x_3)$-plane). Denote by
$g(t)$ a smooth continuation of $h(a,0,t)$ across $t_1$. On the leaves where
$t>t_1$, we now define $h$ (globally) by again requiring (\ref{eq8}), but by
replacing (\ref{eq7}) by $h(a,0,t) =g(t)$. (As above, we use that by
assumption, $\alpha$ is exact on the leaves of $\mathcal{F}$). In order to do
this, we need to know that for $t$ close enough to $t_1$ (and $t>t_1$), these
leaves are distinct from the ones through points $(-a,0,t')$, with $t'>t_0$
and close to $t_0$. Since $L_1$ and $L_2$ have only one hyperbolic point,
namely the origin, in their closure, and no elliptic point (no elliptic point
is contained in the closure of a leaf by the regularity and uniqueness results
in [\ref{K-W}]), their closures intersect only at the origin, and they are
smoothly bounded away from the origin (by an integral curve of the vector field
induced on $\Gamma$ by $\eta/\Gamma$). Consequently, away from the origin,
there is transverse uniformity ``up to the boundary'' (see [\ref{C-N}],
Theorem~3, chapter~III) for $L_1$ and $L_2$. This implies that if $t'$ and $t$
are sufficiently close to $t_0$ and $t_1$, respectively, the above leaves are
indeed distinct. Note that we have used here the hypothesis that $L_1$ and
$L_2$ are distinct. Finally observe that an argument similar to the one used
above to show that $h(a,0,t)$ is smooth up to $t=t_1$ shows that
$x_3$-derivatives of $h$ are bounded near the origin. Since complex
tangential derivatives of $h$ are equal to components of $\alpha$, we get
inductively that near the origin, all derivatives of $h$ are bounded, i.e.\
$h$ is smooth up to the boundary.

We have now extended $h$ into a saturated neighborhood of $P$ (this is a
slight abuse of language since $P$   is not in $\circk$.) The extended
function is $C^\infty$ on this saturated subset of $\circk$ (observe that the
extension of $h$ along the leaves can be obtained via line integrals of
$\alpha$, so that the smoothness of $\alpha$, together with the smoothness of
$h$ near $P$, gives differentiability of $h$ also in directions transverse to
the leaves; compare in this context the discussion of the holonomy maps in
[\ref{C-N}], chapter~IV).

The function $h$ is constructed in a saturated neighborhood of the union of
the hyperbolic leaves (leaves having a hyperbolic point in their closure)
 by repeating the above construction for every hyperbolic point. We
use here the assumption that two leaves meeting at a hyperbolic point have
positive distance from the remaining hyperbolic points. This means that if
$L_1$ and $L_2$ meet at $P_1$, and $L_3$ and $L_4$ meet at $P_2$ then there
are disjoint saturated neighborhoods of $L_1\cup L_2$ and $L_3\cup L_4$ (this
argument is similar to the one used in the paragraph before the last, or see
below, where we show that two leaves whose closures contain no hyperbolic
points can be separated by saturated neighborhoods).  Consequently, this
construction gives a well defined function $h$ that is $C^\infty$ in a
saturated neighborhood of the union of the hyperbolic leaves. 

A similar construction works near the elliptic complex tangencies. Near such a
point, $\circk$ is foliated by a smooth one parameter family of disks, and
there is a (uniformly) transverse  curve emanating from the elliptic point (see
[\ref{K-W}]). Prescribing $h$ to be zero on this curve and then solving
$dh=\alpha$ on each leaf gives the required function $h$ in a saturated subset
of $\circk$ near the elliptic points. Note that the disks close to an elliptic
point do not  contain other complex trangencies in their closure (by the
regularity up to the boundary results in [\ref{K-W}]), so that we choose  the
saturated neighborhoods small enough to be pairwise disjoint and disjoint from
a saturated neighborhood of the union of the hyperbolic leaves. 

For the second step of the construction, we consider the leaf space
$\circk/\mathcal{F}$: \ points in the same leaf are identified, and the
topology is the quotient topology. $\circk/\mathcal{F}$ is in general not
Hausdorff (two hyperbolic leaves whose closure intersects in a hyperbolic
point cannot be separated by saturated neighborhoods and thus give rise
to a pair of points in $\circk/\mathcal{F}$ which cannot be separated). In our
situation, $\circk/\mathcal{F}$ carries a natural structure of a $C^\infty$
one dimensional real (in general non-Hausdorff) oriented manifold ([\ref{H}],
section~3). We use here that the leaves of $\mathcal{F}$ are closed. Also,
note that the first Betti number of $\circk$ is finite (this is assumed in all
of section~3 in [\ref{H}]). The local charts are obtained as follows. Pick a
point $P$ on a leaf $L$ and choose a (short) transversal $C^\infty$ curve
through $P$ that meets each leaf of $\mathcal{F}$ at most once ([\ref{H}],
Proposition~3.1). The projection of this curve gives a neighborhood $U$ of $L$
in $\circk/\mathcal{F}$, and we take  a smooth parameter on the curve (for
instance arclength) to be a local coordinate in $U$. Then the coordinate
transformations are given by holonomy maps and are thus $C^\infty$
([\ref{C-N}], chapter~IV). Because $\mathcal{F}$ is transversely oriented (by
the vector field $T$), $\circk/\mathcal{F}$ is oriented. Note that the
projection $\pi\colon \ \circk\to \circk/\mathcal{F}$ is a $C^\infty$
submersion.

Denote by $L_1,L_2,\ldots, L_m$ the hyperbolic leaves, and set $H := L_1\cup
L_2\cup\cdots\cup L_m$. Since the leaves are closed in $\circk,\circk
\backslash H$ is open. Consider now the quotient $(\circk\backslash
H)/\mathcal{F}$ (which is diffeomorphic to $(\circk/\mathcal{F}) \backslash
\{L_1,\ldots, L_m\}$). We claim that it is Hausdorff. The argument is similar
to one used above. First note that any leaf $L$ in $\circk\backslash H$ is
smooth up to $\Gamma$ and intersects $\Gamma$ in an integral curve of the
vector field induced on $\Gamma$ by   $\eta/\Gamma$. The reason is that
the closure $\bar L$ of $L$ in $K$ contains neither hyperbolic points (by the
definition of $H$) nor elliptic points (no elliptic point is contained in the
closure of a leaf of $\circk$, by the results of [\ref{K-W}]). But near every
other point $P$, $\mathcal{F}$ is smooth up to $\Gamma$ and there is a unique leaf that
intersects $\Gamma$ (from one side) transversally in an integral curve through
$P$ of the field induced by $\eta$. This implies that if $L_1$ and $L_2$ are
two distinct leaves in $\circk\backslash H$, then they (equivalently: \ their
closures) are a positive distance apart. For if not,  their closures (in $K$)
would have to intersect (in points of $\Gamma$). This is impossible by the
preceding discussion. Moreover, for these leaves, there is again transverse
uniformity (see [\ref{C-N}], Theorem~3, chapter III) ``up to the boundary''.
Consequently, they can be separated by {\em saturated\/} neighborhoods. This
proves that $(\circk\backslash H)/\mathcal{F}$ is Hausdorff. With the
differential structure induced as in the case of the full leaf space
(equivalently: \ the restriction of the differential structure on
$\circk/\mathcal{F}$ to $\circk/\mathcal{F} \backslash\{L_1,\ldots, L_m\}$),
$(\circk\backslash H)/\mathcal{F}$ is thus a bona fide real 1-dimensional
manifold. By the classification of such manifolds, its connected components
are diffeomorphic either to $\mathbb{R}$ or to $S^1$.

We now construct $h$ on the components of $\circk\backslash H$ (which
correspond to components of $(\circk\backslash H)/\mathcal{F}$) as follows. Let
$U$ be such a component. Then by the above, $\pi(U)$ is diffeomorphic to
$\mathbb{R}$ or to $S^1$. Denote by $\gamma$ the one-form on $\pi(U)$
corresponding to $dx$ in the case of $\mathbb{R}$, and $d\theta$ in the case
of $S^1$. Note that $\gamma$ is closed and nonvanishing. The form $\omega :=
\pi^*\gamma$ is then a nonsingular ($\pi$ is a submersion) one-form on $U$ that
satisfies
\begin{equation}\label{eq9}
d\omega = d\pi^*\gamma = \pi^*d\gamma = 0,
\end{equation}
and
\begin{equation}\label{eq10}
\omega/_L = 0,
\end{equation}
where $L$ is any leaf of $\mathcal{F}/_U$. Thus $\omega$ is of the form $
e^g\eta$ for some $C^\infty$-function $g$ (where $e^g = \omega(T)$). Therefore,
\begin{equation}\label{eq11}
0 = d\omega = e^g(dg\wedge \eta + d\eta) = e^g(dg\wedge \eta + \alpha \wedge
\eta).
\end{equation}
Here we have used that on $K, d\eta = \alpha\wedge\eta$ (see the discussion at
the end of section~2). Consequently,
\begin{equation}\label{eq12}
(dg+\alpha) \wedge \eta = 0,
\end{equation}
which is equivalent to
\begin{equation}\label{eq13}
dg/_L = -\alpha/_L
\end{equation}
on each leaf of $\mathcal{F}/_U$. (This computation will also show the
equivalence of (i) and (ii) in Proposition~\ref{pro2} below.)

In the last step, we patch together the functions $h$ and $g$ constructed so
far. By working in the leaf space (equivalently: \ on (short) integral curves
of $T$), one constructs a smooth function $\psi$ which is constant on the
leaves of $\mathcal{F}$, is identically one in a (saturated) neighborhood of
$L_1\cup\cdots\cup L_m$, is supported in a small (saturated) neighborhood of
this set, is identically one near elliptic points and is supported close to
these points. Then the function $h_1 := \psi h + (1-\psi)g$ satisfies, for a
leaf $L$ of $\mathcal{F}$, 
\begin{align}
dh_1/_L = d(\psi h + (1-\psi)g/_L &= (\psi dh + (1-\psi)dg)/_L\nonumber\\
\label{eq14}
&= \psi dh/_L + (1-\psi) dg/_L\\
&= \psi(-\alpha)/_L + (1-\psi) (-\alpha)/_L = (-\alpha)/_L.\nonumber
\end{align}
That is, $h_1$ satisfies (\ref{eq6}).

It remains to check smoothness up to the boundary. As indicated at the
beginning of the proof of Theorem~\ref{thm1}, we will verify that derivatives
of all orders are bounded. We have already  noted that by our construction in
local coordinates near a hyperbolic point of $\Gamma$, this is the case. The
essence of the argument is the same near the other points of $\Gamma$. Let
$P$ be a point of $\Gamma$ which is not a complex tangency. Near $P$, the Levi
foliation is smooth up to $\Gamma$, the boundaries of the leaves being
integral curves of the vector field on $\Gamma$ induced by $\eta/_\Gamma$.
Denote by $L_p$ the leaf ``through'' $P$. For $Q\in L_p$, sufficiently close
to $P$, let $\sigma(t)$, $-\vp < t<\vp$, be a short integral curve of $T$
through $Q$ (say $\sigma(0) = Q$). For $z\in \circk$ near $P$, denote by $L_z$
the leaf of $\mathcal{F}$ through $z$, and by $t_z$ the unique $t$ such that
$\sigma(t_z)\in L_z$. Writing $h_1(z)$ as $h_1(\sigma(t_z)$) plus a line
integral of $\alpha$ in $L_z$ shows, in view of the smoothness of $\alpha$ and
of $\mathcal{F}$  up to $\Gamma$, that derivatives of all orders of $h_1$ are
bounded near $P$. Analogous arguments work near an elliptpic point. This
completes the proof of Theorem~\ref{thm1}. \medskip

\begin{rem}\label{rem6}
Near an elliptic point, one can take advantage of the
fact that conditions (\ref{eq1}) and (\ref{eq2}) need only be satisfied to
within $\vp$. Specifically, if $P$ is an elliptic point, define $h_\vp$ near
$P$ through $e^{h_\vp} L_n = L_n(P)$. Then, close enough to $P$, $|\arg
e^{h_\vp} L_n(\rho)|<\vp$, and $\partial\rho([e^{h_\vp}L_n,\partial/\partial
\bar z_j])=0$ (the latter holds since $e^{h_\vp} L_n$ has holomorphic (namely
constant) coefficients). 
\end{rem}\bigskip

\centerline{4. MORE FOLIATION THEORY}

The existence of the function $h$ constructed above also entails some
properties of interest from the point of view of foliation theory, both of the
Levi foliation and of the flow generated by the transverse field $T$. This
leads to to a homological condition that is necessary and sufficient for the
construction from section~3, and thus for the vector fields from section~2,
with $\vp=0$ on $K$, that uses Sullivan's theory of foliation currents
([\ref{Su1}], [\ref{Su2}], [\ref{C-C}]). Although this does not seem to
simplify the proof of, for example, Theorem~\ref{thm1}, it is of independent
interest, and we expect this point of view to be useful in the more general
situation where the set $K$ of infinite type points is not the closure of a
nicely foliated open subset of the boundary.

We keep the notation from section~2. For the moment, we let $M$ be a
$C^\infty$-smooth (open) Levi flat hypersurface in $\mathbb{C}^n$.
$\mathcal{F}$ again denotes the Levi foliation on $M$ (in this case by complex
manifolds of dimension $n-1$).

A one form $\omega$ is said to define the foliation $\mathcal{F}$ if for each
$P\in M$, the tangent space to the leaf through $P$ is equal to the kernel of
$\omega$ at $P$. So the form $\eta$ defined in section~2 always defines
$\mathcal{F}$.

A flow is said to be geodesible if there exists a Riemannian metric (on the
underlying manifold) such that the orbits become geodesics (see e.g.\
[\ref{T}], chapter~6). We consider the flow generated by $T$.

\begin{pro}\label{pro2}
The following are equivalent:
\begin{itemize}
\item[(i)] there exists $h\in C^\infty(M)$ such that
$$dh/_L + \alpha/_L = 0,\qquad \forall L\in \mathcal{F}$$
\item[(ii)] there exists a smooth $(C^\infty)$ {\em closed\/} one form
$\omega$ on $M$ that defines $\mathcal{F}$
\item[(iii)] the flow generated by $T$ is geodesible by a metric that is
obtained from the Euclidean metric on $M$ by scaling in the $T$-direction.
\end{itemize}
\end{pro}

Note that (ii) and (iii) are manifestly independent of a choice of defining
function, while the definition of $\alpha$ involves such a choice. However, we
have pointed out already in section~2 that one can check by direct computation
that (i) is actually independent of this choice.

\begin{proof}
As noted in section~3 the equivalence of (i) and (ii) follows form the
computation (\ref{eq11})--(\ref{eq13}) there. Namely, since $\omega$
defines $\mathcal{F}$, it is nonsingular, so $\omega=e^h\eta$ for some real
valued function $h$ (possibly after replacing $\omega$ by $-\omega$). Then 
\begin{equation}\label{eq15}
\begin{split}
d\omega = d(e^h\eta) &= e^h(dh\wedge \eta + d\eta)\\
&= e^h(dh+\alpha)\wedge \eta.\\
\end{split}
\end{equation}
Therefore, $d\omega=0 \Leftrightarrow (dh+\alpha) \wedge
\eta=0\Leftrightarrow$ (i) holds. In the last equality in (\ref{eq15}), we have
again used that $d\eta = \alpha\wedge \eta$ (see section~2).

The equivalence of (ii) and (iii) is contained in Proposition~6.7 of
[\ref{T}], where various characterizations are given for a flow to be
geodesible. But the arguments show that (ii) is equivalent to the flow of $T$
being geodesible by a metric as in (iii). This completes the proof of
Proposition~\ref{pro2}.\medskip

Sullivan ([\ref{Su2}]) characterizes geodesible flows by a homological
condition. His arguments are easily modified to obtain a characterization of
the conditions in Proposition~2 (actually, the
versions on the compact set $K$, rather than on an open hypersurface, as in
Proposition~2). We return to the analogue of the setup in section~3 for a
domain in $\mathbb{C}^n$. Denote by $\mathcal{T}$ the one dimensional
foliation whose leaves are the orbits of $T$. (This is actually a foliation of
all of $b\Omega$). We briefly recall the notions we need from the theory of
foliation currents; for details, see [\ref{Su1}] or [\ref{C-C}], chapter 10.
Denote by $\mathcal{D}_1(K)$ the space of one forms that are $C^\infty$ on
$K$, with the usual Fr\'echet space topology. Its strong dual
$\mathcal{D}_1(K)'$ comprises the one currents on $K$. $\mathcal{D}_1(K)$ and
$\mathcal{D}_1(K)'$ are mutually dual (i.e.\  $\mathcal{D}_1(K)$ is
reflexive). The closed convex cone $C_\tau$ of foliation currents for
$\mathcal{T}/_K$ is the closure, in $\mathcal{D}_1(K)'$, of the set of all
finite linear combinations with positive coefficients of currents given by
pairing at $x$ with $T(x)$, for $x\in K$ (so-called Dirac currents). The
crux of the matter is that $C_\tau$ has a compact base: \ $C_\tau\cap
\eta^{-1}(0) = \{0\}$ and $C_\tau \cap \eta^{-1}(1)$ is compact (where $\eta$
is viewed as an element of $D_1(K)''$); this is analogous to Lemma~10.2.3 in
[\ref{C-C}]. We say that a $p$-chain is tangential to $\mathcal{F}$ if it is a
sum of $p$-chains, each of which is a $p$-chain in a leaf of
$\mathcal{F}$. \end{proof}

\begin{thm}\label{thm3}
There exists $h\in \mathcal{D}(K)$ satisfying $dh/_L + \alpha/_L = 0$ for each
leaf $L$ of $\mathcal{F}$  if and only if no nontrivial
foliation current for $\mathcal{T}/_K$ can be arbitrarily
well approximated in $\mathcal{D}_1(K)'$ (i.e.\ uuniformly on bounded subsets
of $\mathcal{D}_1(K)$) by sums of the form (smooth 1-chain tangential to
$\mathcal{F}$) plus (boundary of smooth 2-chain in $\circk$).
\end{thm}
\medskip

The proof of Theorem~\ref{thm3} is completely analogous to the proof of
(iii) in the theorem in [\ref{Su2}]. Namely, the second statement in
Theorem~3 means that the closed linear subspace of $\mathcal{D}_1(K)'$
generated by boundaries of smooth 2-chains in $\circk$ and smooth 1-chains
tangential to $\mathcal{F}$ intersects the cone $C_\tau$ only in $0\in
\mathcal{D}_1(K)'$. Since $C_\tau$ has a compact base (see above), a suitable
version of the Hahn-Banach theorem (see e.g.\ [\ref{C-C}], Theorem~10.2.5 and
the remark immediately following it) yields an element $\omega\in
\mathcal{D}_1(K)'' = \mathcal{D}_1(K)$ that vanishes on he boundaries of
smooth 2-chains and on smooth 1-chains tangential to $\mathcal{F}$, while it
is strictly positive on $C_\tau\backslash\{0\}$. These properties imply,
respectively, that $\omega$ is closed, vanishes at all points on vectors
tangent to the leaf of $\mathcal{F}$ through the point, and is strictly
positive on $T$ at all points of $K$. Consequently, $\omega = e^h\eta$, with
$h\in \mathcal{D}(K)$, and $h$ has the required property, see (the proof of)
Proposition~\ref{pro2}. 

On the other hand, if there exists $h$ as in the
theorem, then $\omega = e^h\eta$ is closed (again by the proof of
Proposition~\ref{pro2}), and, if $\varphi\in C_\tau$, then $\varphi(\omega)>0$,
so that $\varphi$ cannot be approximated as indicated in the theorem. This
concludes the proof of Theorem~\ref{thm3}.

\begin{rem}\label{rem7}
As an illustration, consider the following.
We know from section~3 that if there exists a function $h$ as in
Theorem~\ref{thm3}, then the Levi foliation must be without holonomy. It
is easily seen directly that the homological condition in Theorem~\ref{thm3}
implies that the holonomy of $\mathcal{F}$ is trivial. Consider a point $x_0$
in a leaf $L$ and a simple loop $\sigma$ in $L$, starting and ending at $x_0$.
Lift $\sigma$ to a nearby leaf as in the construction of the holonomy maps.
This yields a path $\tilde\sigma$. If $\tilde\sigma$ is not a loop, we have a
two dimensional surface whose boundary consists of $\sigma$ plus $\tilde\sigma$
(both of which are tangential to the Levi foliation) plus a piece of an
integral curve of $T$ (i.e.\ a nontrivial foliation current for $\mathcal{T}$).
This is impossible by the homological condition. Since $\sigma$ was arbitrary,
$L$ has trivial holonomy.
\end{rem}

\newpage

\centerline{REFERENCES}

\begin{enumerate}
\item\label{Ba1}
David Barrett: \ Complex analytic realization of Reeb's foliation of $S^3$,
Math.\ Z. {\bf 203} (1990), 355--361.

\item\label{Ba-2}
David Barrett:\ Behavior of the Bergman projection on the Diederich-Fornaess
worm, Acta Math. {\bf 168}, 1--2 (1992), 1--10.

\item\label{B-F1} 
Eric Bedford and John Erik Fornaess:\ Domains with pseudoconvex neighborhood
systems, Invent.\ Math. {\bf 47} (1978), 1--27.

\item\label{B-G}
Eric Bedford and Bernard Gaveau:\ Envelopes of holomorphy of certain 2-spheres
in $\mathbb{C}^2$, American J.\ of Math. {\bf 105} (1983), 975--1009.

\item\label{B} 
Errett Bishop:\ Differentiable manifolds in complex Euclidean space, Duke
Math.\ J. {\bf 32} (1965), 1--21.

\item\label{B-S3}
Harold P.\ Boas and Emil J.\ Straube: \ Sobolev estimates for the
$\bar\partial$-Neumann operator on domains in $\mathbb{C}^n$ admitting a
defining function that is plurisubharmonic on the boundary, Math.\ Z. {\bf
206} (1991), 81--88.

\item\label{B-S1}
Harold P.\ Boas and Emil J.\ Straube: \ De Rham cohomology of manifolds
containing the points of infinite type, and Sobolev estimates for the
$\bar\partial$-Neumann problem, J.\ Geometric Analysis {\bf 3}, Nr.~3 (1993),
225--235.

\item\label{B-S2}
Harold P.\ Boas and Emil J.\ Straube: \ Global regularity of the
$\bar\partial$-Neumann problem: \ a survey of the $\mathcal{L}^2$-Sobolev
theory, Several Complex Variables, M.\ Schneider and Y.-T.\ Siu editors,
Mathematical Sciences Research Institute Publications {\bf 37}, 79--111,
Cambridge Univ.\ Press, 1999.

\item\label{C-N}
C\'esar Camacho and Alcides Lins Neto:\ Geometric Theory of Foliations,
Birkh\"auser, 1985.

\item\label{C-C} 
Alberto Candel and Lawrence Conlon:\ Foliations I, Graduate Studies in
Mathematics, vol.~23, American Math. Society, 2000.

\item\label{DA} 
John P.\ D'Angelo:\ Several Complex Variables and the Geometry of
Real Hypersurfaces, CRC Press, 1993.

\item\label{F-K}
G.B.\ Folland and J.J.\ Kohn: \ The Neumann Problem for the Cauchy-Riemann
Complex, Annals of Math.\ Studies {\bf 75}, Princeton Univ.\ Press, 1972.

\item\label{F}
Franc Forstneri\v c:\ Complex tangents of real surfaces in complex surfaces,
Duke Math.\ J. {\bf 67} (1992), 353--376.

\item\label{H}
Andr\'e H\"afliger:\ Vari\'et\'es feuillet\'ees, Ann.\ Scuola Norm.\ Sup.\ Pisa
{\bf 16} (1962), 367--397.

\item\label{K-W}
Carlos E.\ Kenig and Sidney M.\ Webster:\ The local hull of holomorphy of a
surface in the space of two complex variables, Invent.\ Math. {\bf 67} (1982),
1--21.

\item\label{S-S}
Emil J.\ Straube and Marcel K.\ Sucheston: \ Plurisubharmonic defining
functions, good vector fields, and exactness of a certain one form, in
preparation.

\item\label{Su1}
Dennis Sullivan: \ Cycles for the dynamical study of foliated manifolds and
complex manifolds, Invent.\ Math. {\bf 36} (1976), 225--255.

\item\label{Su2}
Dennis Sullivan: \ A foliation of geodesics is characterized by having no
``tangent homologies'', J.\ Pure \& Applied Algebra {\bf 13} (1978), 101--104.

\item\label{T} 
Philippe Tondeur:\ Geometry of Foliations, Birkh\"auser, 1997.
\end{enumerate}

\vspace{.5in}

\baselineskip = 12pt

\hspace{2.5in} Department of Mathematics

\hspace{2.5in} Texas A\&M University

\hspace{2.5in} College Station, TX \ 77843-3368

\hspace{2.5in} straube@math.tamu.edu

\end{document}